\documentclass[a4paper,12pt]{article}
\usepackage{listings}
\usepackage{graphicx}
\usepackage{subcaption}
\usepackage{color}
\usepackage[ruled,vlined]{algorithm2e}
\usepackage{amsmath, amsthm}
\usepackage{amssymb,amsfonts}
\usepackage{enumerate}
\usepackage{csquotes}
\usepackage{fancyhdr}
\usepackage{graphicx}
\usepackage{hyperref}
\usepackage{authblk}

\usepackage{amssymb}

\theoremstyle{definition}
\newtheorem{definition}{Definition}[section]
\newtheorem{example}[definition]{Example}

\newtheorem{remark}[definition]{Remark}

\theoremstyle{plain}
\newtheorem{theorem}[definition]{Theorem}
\newtheorem{lemma}[definition]{Lemma}

\newtheorem{corollary}[definition]{Corollary}
\numberwithin{equation}{section}

\usepackage{color}   
\usepackage{hyperref}
\usepackage{url}
\begin{document}
\date{}
\author[1]{Mr. Wajid M. Shaikh}
\author[2]{Dr. Rupali S. Jain}
\author[3]{Dr. B. Surendranath Reddy}
\affil[1]{Research Scholar, School of Mathematical Sciences, S.R.T.M. University, Nanded}
\affil[2]{Associate Professor, School of Mathematical Sciences, S.R.T.M. University, Nanded}
\affil[3]{Assistant Professor, School of Mathematical Sciences, S.R.T.M. University, Nanded}

\title{
     Construction of Linear Codes from the Unit Graph $G(\mathbb{Z}_{n}\oplus \mathbb{Z}_{m})$}
\maketitle
\begin{abstract}
In this paper, we develop the python code for generating unit graph $G(\mathbb{Z}_{n}\oplus\mathbb{Z}_{m})$, for any integers $m\ \& \ n$.  For any prime $r$, we construct $r$-ary linear codes from the incidence matrix of the unit graph $G(\mathbb{Z}_{n}\oplus\mathbb{Z}_{m})$, where $n \ \& \ m$ are either power of prime or product of power of primes. We also prove the minimum distance of dual of the constructed codes as  either 3 or 4. Finally, we state conjectures two on linear codes constructed from the unit graph $G(\mathbb{Z}_{n}\oplus \mathbb{Z}_{m})$, for any integer $m\ \&  \ n$.
\end{abstract}

\section{Introduction}
\par The specific applications of linear codes over finite fields in computer and communication systems, data storage devices, and consumer electronics have been thoroughly examined by researchers. So construction of linear  codes from different types of functions and graphs have been extensively studied. In this paper, we construct linear codes from unit graphs $G(\mathbb{Z}_{n}\oplus \mathbb{Z}_{m})$.
\par In 1990, the unit graph over $\mathbb{Z}_{n}$, introduced by Grimaldi, R.P. \cite{14}, is the simple graph, where $\displaystyle x,y\in \mathbb{Z}_{n}$ are adjacent if and only if $x+y$ is a unit in $\mathbb{Z}_{n}$.
 \par Construction of linear codes from incidence matrices of line graphs and Hamming graphs was studied by Fish, W., Key, J. D., $\&$ Mwambene, E.\cite{7}. Also, Key, J. D., $\&$ Rodrigues, B. G. \cite{8} utilized lattice graphs to construct codes and investigated their decoding techniques through permutation decoding. In 2013, Dankelmann, P., Key, J. D., $\&$ Rodrigues, B. G. \cite{3} provided a generalization of the relationship between the parameters of connected graphs and the linear codes derived from the incidence matrices of these graphs.  Many researchers have extensive research on constructing linear codes using adjacency matrices of specific graphs \cite{11,12,13}.
\par Recently, Annamalai, N., $\&$ Durairajan, C.\cite{1},\cite{21}, constructed linear codes from the incidence matrices of unit graphs and zero divisor graphs. These results were generalized by Jain, R. S., Reddy, B. S. \& Shaikh, W. M. \cite{20}.  In this paper, we extended their work by developing python code for generating unit graph $G(\mathbb{Z}_{n}\oplus \mathbb{Z}_{m})$, for any integer $n,m$. Also, we construct linear codes from the  incidence matrices of unit graph $\displaystyle G(\mathbb{Z}_{n}\oplus\mathbb{Z}_{m})$ and find the parameters  for  their dual codes.  Finally we conclude by stating two conjectures.
\section{Preliminaries}
In this section, we recall definitions and results related to  graphs and linear codes. Let $\mathbb{Z}_{n}\oplus \mathbb{Z}_{m}$ denote the ring under the component wise modulo addition and modulo multiplication. Here, we denote units and non-units of $\displaystyle \mathbb{Z}_{n}\oplus \mathbb{Z}_{m}$ by $\displaystyle U(\mathbb{Z}_{n}\oplus\mathbb{Z}_{m})$ and $\displaystyle N_{U}(\mathbb{Z}_{n}\oplus \mathbb{Z}_{m})$ respectively. Note that, $\bar{x}=(x_{1},x_{2})$ is unit in $\mathbb{Z}_{n}\oplus \mathbb{Z}_{m}$ if and only if both $x_{1}$ and $x_{2}$ are units in $\mathbb{Z}_{n}$ and $
\mathbb{Z}_{m}$ respectively.
\begin{definition}\cite{18}[Linear Code]
Let $\displaystyle \mathbb{F}_{r}$ represents the finite field with $r$ elements. A  linear code $C_{r}$ of length $n$ is a subspace of $\displaystyle \mathbb{F}^{n}_{r}$ and it is called $r$-ary linear code. Dimension of linear code $C_{r}$ is the dimension of $C_{r}$ as a vector space over field  $\displaystyle \mathbb{F}_{r}$ and is denoted by $\displaystyle \text{dim}(C_{r})$.
\end{definition}
\begin{definition}\cite{18}[Dual of code]
Let $C_{r}$ be a linear code of length $n$ over $\displaystyle \mathbb{F}_{r}$. Then dual of code $C_{r}$ is the orthogonal compliment of the subspace $C_{r}$ in $\displaystyle \mathbb{F}^{n}_{r}$ and is denoted by $\displaystyle C^{\perp}_{r}$.
\end{definition}
\begin{theorem}\cite{18}\label{thm8}
Let $C_{r}$ be a $r$-ary code of length $n$ over a field $\displaystyle \mathbb{F}_{r}$. Then $\displaystyle C^{\perp}_{r}$ is a linear code of length $n$ and $\displaystyle \text{dim}(C_{r}^{\perp})=n-\text{dim}(C_{r})$.
\end{theorem}

\begin{definition}\cite{18}[Minimum Hamming distance]
Let $C_{r}$ be a linear code. The minimum Hamming distance of code $C_{r}$, denoted by $d(C_{r})$, is defined as
$\displaystyle d(C_{r})=\text{min}\{d_{C}(x,y) \ | \ x,y\in C_{r} \ \& \ x\neq y\}$,\\
 where $d_{C}(x,y)=\text{Number of places where $x$ and $y$ differs}$.
\end{definition}

\begin{remark}\cite{18}
A $r$-ary linear code $C_{r}$ of length $n$, dimension $k$ and minimum distance $d$ is called $[n,k,d]_{r}$ linear code.
\end{remark}
\begin{definition}\cite{18}
A generator matrix of linear code $C_{r}$ is a matrix $H$ whose rows form a basis for $C_{r}$ and a generator matrix $\displaystyle H^{\perp}$ of linear code $C^{\perp}_{r}$ is called parity-check matrix of $C_{r}$.
\end{definition}
\begin{definition}\cite{19}
The distance between two vertices $x$ and $y$, denoted by $d(x,y)$, is the length of a shortest path from $x$ to $y$. The diameter of a graph $G$ is denoted by $\text{diam}(G)$, is the maximum distance between any two vertices in $G$. i.e.
$\displaystyle \text{diam}(G)=\text{Max}\{d(x,y) \ | \ x,y\in V\}$.
\end{definition}

\begin{definition}\cite{19}
Let $G$ be a simple graph. The edge connectivity of $G$, denoted by $\lambda(G)$, is the smallest number of edges in $G$ whose deletion from $G$ either leaves a disconnected graph or an empty graph.
\end{definition}
\begin{definition}\cite{2}
	Let $\mathcal{R}$ be a ring with nonzero identity. The unit graph of $\mathcal{R}$, denoted by $G(\mathcal{R})$, is a graph with vertex set as $
	\mathcal{R}$ and two distinct vertices $x$ and $y$ are adjacent if and only if $x + y$ is a unit of $\mathcal{R}$.
\end{definition}
\begin{theorem}\label{thm3}\cite{4}
Let $\displaystyle G$ be a connected graph with vertex set $V$.\\ If $\text{diam}(G)\leq 2$ then edge connectivity of $G$ is  $\lambda(G)=\delta(G)$.
\end{theorem}

\begin{theorem}\label{thm4}\cite{5}
Let $\displaystyle G$ be a connected bipartite graph, if $\displaystyle \text{diam}(G)\leq 3$ then edge connectivity of $G$ is $\displaystyle \lambda(G)=\delta(G)$.
\end{theorem}

\begin{theorem}\label{thm2}\cite{3}
Let $G$ be a connected graph and let $H$ be a $\displaystyle |V|\times |E|$ incidence matrix for $G$. Then binary code generated by $H$ is\\ $\displaystyle C_{2}(H)=\displaystyle [|E|,|V|-1,\lambda(G)]_{2}$.
\end{theorem}

\begin{theorem}\cite{3}\label{thm5}
  Let $G$ be a connected bipartite graph and let $H$ be a $\displaystyle |V|\times |E|$ incidence matrix for $G$, and $r$ be an odd prime. Then $r$-ary code generated by $H$ is $\displaystyle C_{r}=[|E|,|V|-1,\lambda(G)]_{r}$.
\end{theorem}

\begin{theorem}\cite{3}\label{thm6}
Let $G$ be a connected graph with girth $g_{r}(G)$ and even girth $g_{r}(G)_{e}$. Let $H$ be an incidence matrix for $G$, $C=C_{r}(H)$, where $r$ is any prime, and $\displaystyle d^{\perp}$ is the minimum distance of $\displaystyle C_{r}^{\perp}$. If $q=2$ or $g_{r}(G)$ is even then $\displaystyle d^{\perp}=g_{r}(G)$.
\end{theorem}
\begin{theorem}\cite{20}\label{thm9}
Let $\displaystyle G(\mathbb{Z}_{2^{m}p^{n}})$ be a unit graph, where $\displaystyle p$ is an odd prime. Then
\begin{enumerate}
	\item $\displaystyle |V|=2^{m}p^{n}$ and $\displaystyle |E|=2^{m-1}p^{n}\phi(2^{m}p^{n})$.
	\item  $\displaystyle \lambda(G(\mathbb{Z}_{2^{m}p^{n}}))=\delta(G(\mathbb{Z}_{2^{m}p^{n}}))=\phi(2^{m}p^{n})$.
\end{enumerate}
\end{theorem}
\begin{corollary}\cite{20}\label{thm10}
	Let $\displaystyle G(\mathbb{Z}_{2^{m}p^{n}})$ be a unit graph, where $p$ is an odd prime.
	\begin{enumerate}
		\item If $\displaystyle 2^{m}p^{n}\neq 6$ then $\displaystyle g_{r}(G(\mathbb{Z}_{2^{m}p^{n}}))=4$
		\item If $\displaystyle 2^{m}p^{n}= 6$ then $\displaystyle g_{r}(G(\mathbb{Z}_{2^{m}p^{n}}))=6$
	\end{enumerate}
\end{corollary}
\section{$G(\mathbb{Z}_{n}\oplus \mathbb{Z}_{m})$ and Python Code}
In this section, to visualize the structure of unit graph $G(\mathbb{Z}_{n}\oplus \mathbb{Z}_{m})$, we have programmed the python code for generating unit graph $G(\mathbb{Z}_{n}\oplus \mathbb{Z}_{m})$ for any values of $m$ and $n$. We also, discuss some examples of unit graphs generated from the python code. Lastly, we give the theorem on number of edges of unit $G(\mathbb{Z}_{n}\oplus \mathbb{Z}_{m})$ and condition for the graph $G(\mathbb{Z}_{n}\oplus \mathbb{Z}_{m})$ to be bipartite.

To start with, let us consider an example of unit graph $G(\mathbb{Z}_{5}\oplus \mathbb{Z}_{5})$.  This graph has $25$  vertices and $192$ edges. The number of vertices and edges in $G(\mathbb{Z}_{n}\oplus\mathbb{Z}_m)$ increases rapidly even with small increase in $n \ \&  \ m$.
\par
\begin{figure}[h]
	\centering
	\includegraphics[width=0.5\textwidth]{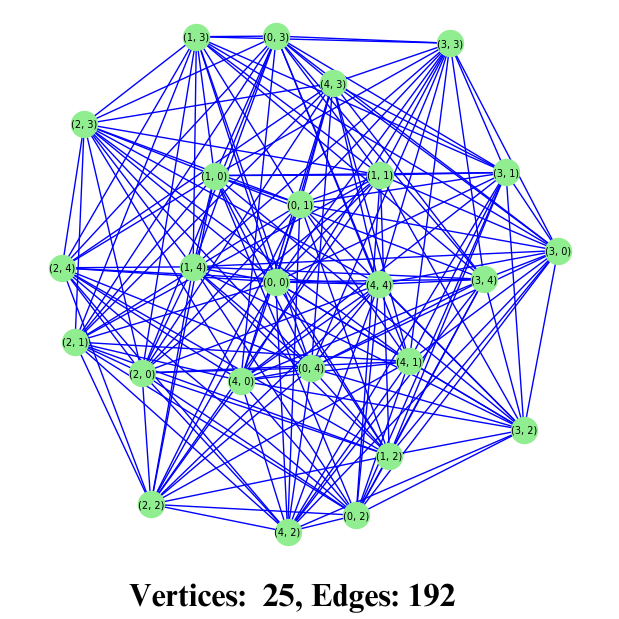} 
	\caption{Unit graph of $G(\mathbb{Z}_{5}\oplus \mathbb{Z}_{5})$}
	\label{fig1} 
\end{figure}

 It is quite difficult for visualizing the structure of unit graph. For this, we develop the python code, which take input as integer values $n$ and $m$ and give output as graph of unit graph, for example $n=5 \ \& \ m=5$ give the unit graph $G(\mathbb{Z}_{n}\oplus \mathbb{Z}_{m})$ as shown in Fig \ref{fig1}.

\begin{lstlisting}
# Importing necessary libraries
import numpy as np
import networkx as nx
import matplotlib.pyplot as plt
# Define modolu multiplication
def dot_product(A,B):
C= ((A[0]*B[0])%n,(A[1]*B[1])%m)
return C
# Define modolu addition
def additionlist(A,B):
C=((A[0]+B[0])%n,(A[1]+B[1])%m)
return C
#For input values of m and n
n=int(input("Please Enter Value of n: "))
m=int(input("Please Enter Value of m: "))
# Generating Z_m, Z_n and Z_m direct sum Z_n
z_m=[]
for i in range(m):
z_m.append(i)
z_n=[]
for i in range(n):
z_n.append(i)

z_nz_m=[]
for i in z_n:
for j in z_m:
z_nz_m.append((i,j))

# Finding non units
Non_Unit=[]
z_nz_m.remove((0,0))
for i in z_nz_m :
for j in z_nz_m:
if (dot_product(i,j)==(0,0) and i not in Non_Unit):
Non_Unit.append(i)

z_nz_m.append((0,0))
Non_Unit.append((0,0))
# Defining edges
edges=[]
for i in z_nz_m:
for j in z_nz_m:
if i!=j and additionlist(i,j) not in Non_Unit:
edges.append((i,j))
# Generating Grpah
graph = nx.Graph()

for name in z_nz_m:
graph.add_node(name, label=name)
graph.add_edges_from(edges)
labels = nx.get_node_attributes(graph, "label")
plt.figure(figsize=(6, 6))
nx.draw(graph, with_labels=True, labels=labels, node_color="lightgreen", node_size=350, edge_color="blue", font_size=7)
plt.show()
\end{lstlisting}
\par In addition to above, the following python code can check the given graph is bipartite or not, as shown in Fig \ref{fig2}.
\begin{lstlisting}
	bipartite = nx.is_bipartite(graph)
	
	if bipartite:
	print("The graph is bipartite.")
	else:
	print("The graph is not bipartite.")
	top_nodes, bottom_nodes = nx.bipartite.sets(graph)
	plt.figure(figsize=(8, 6))
	pos = nx.spring_layout(graph)
	nx.draw_networkx(bipartite_graph, pos=pos, with_labels=True, node_color=["lightgreen" if node in top_set else "pink" for node in bipartite_graph.nodes()], node_size=350, edge_color='blue', font_size=7)
	
	plt.axis("off")
	plt.show()
\end{lstlisting}
\par Also, the following python code gives the, incidence matrix of the unit graph $G(\mathbb{Z}_{n}\oplus \mathbb{Z}_{m})$.
\begin{lstlisting}
np.set_printoptions(threshold=np.inf)
# Get the adjacency matrix of the graph
adjacency_matrix = nx.adjacency_matrix(graph).todense()

# Create a list of edges
edges = list(graph.edges())

# Create a list of vertices
vertices = list(graph.nodes())

# Create an empty incidence matrix
incidence_matrix = np.zeros((len(edges), len(vertices)))

# Fill in the incidence matrix
for i, edge in enumerate(edges):
for j, vertex in enumerate(vertices):
if edge[0] == vertex or edge[1] == vertex:
incidence_matrix[i][j] = 1

# Print the incidence matrix
print(incidence_matrix)
\end{lstlisting}
\begin{example}
Consider the following examples of unit graphs
\begin{enumerate}
\item Unit graph $G(\mathbb{Z}_{4}\oplus \mathbb{Z}_{5})$ is a connected bipartite graph as shown in Figure \ref{fig2}.
\item Unit graph $G(\mathbb{Z}_{6}\oplus \mathbb{Z}_{4})$ is a disconnected graph as shown in Figure \ref{fig3}.
\end{enumerate}
\begin{figure}[h]
	\centering
	\begin{subfigure}[h]{0.45\textwidth}
		\centering
		\includegraphics[width=\linewidth]{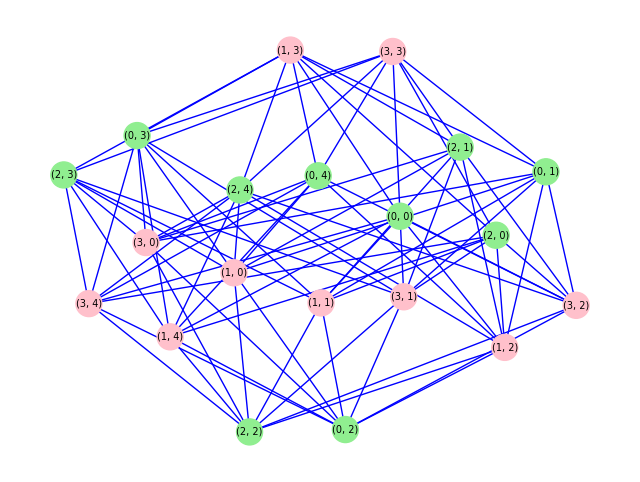}
		\caption{Unit graph of $G(\mathbb{Z}_{4}\oplus \mathbb{Z}_{5})$}
		\label{fig2}
	\end{subfigure}
	\hfill
	\begin{subfigure}[h]{0.45\textwidth}
		\centering
		\includegraphics[width=\linewidth]{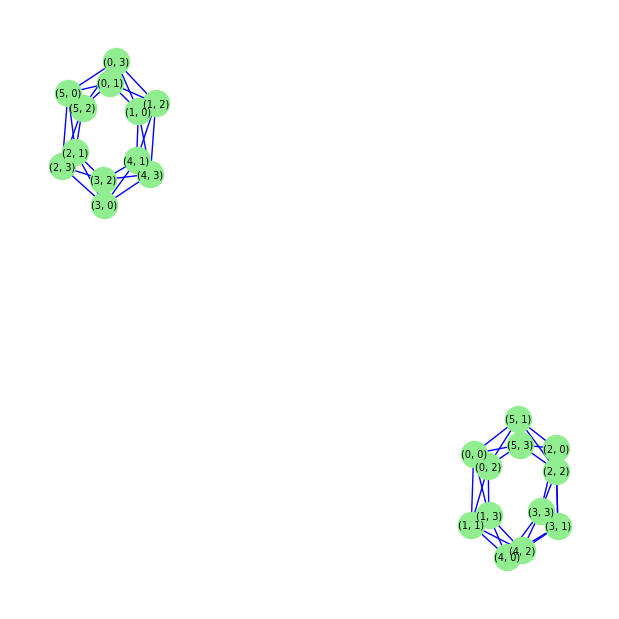}
		\caption{Unit graph of $G(\mathbb{Z}_{6}\oplus \mathbb{Z}_{4})$}
		\label{fig3}
	\end{subfigure}
\end{figure}
\end{example}

\begin{theorem}\label{thm41}
	Let $G(\mathbb{Z}_{n}\oplus \mathbb{Z}_{m})$ be a unit graph, where $m$ and $n$ are any positive integer. Then
	$$\displaystyle|E|=\begin{cases}
		\frac{(mn-1)\phi(m)\phi(n)}{2}, & \text{if both $m$ and $n$ are odd}\\
		\frac{mn\phi(m)\phi(n)}{2}, & \text{otherwise}
	\end{cases}$$
		\begin{proof}
			 First consider, both $m$ and $n$ are odd. If $\bar{x}\in \displaystyle \mathbb{Z}_{n} \oplus \mathbb{Z}_{m} $, then $\bar{y}\in \displaystyle \mathbb{Z}_{n} \oplus \mathbb{Z}_{m}$ is adjacent to $\bar{x}$ if and only if $\bar{x}+\bar{y}\in U(\mathbb{Z}_{n} \oplus \mathbb{Z}_{m})$, if and only if $\bar{y}\in U(\mathbb{Z}_{n} \oplus \mathbb{Z}_{m})-\bar{x} $.
			Note that, if $\bar{x}\in U(\mathbb{Z}_{n} \oplus \mathbb{Z}_{m}) $, then $2\cdot \bar{x} \in U(\mathbb{Z}_{n} \oplus \mathbb{Z}_{m}) $ which implies
			$\bar{x}\in U(\mathbb{Z}_{n} \oplus \mathbb{Z}_{m})-\bar{x} $. Hence, $\text{deg}(\bar{x})=|U(\mathbb{Z}_{n} \oplus \mathbb{Z}_{m}|-1$.\\
			If $\bar{x}\notin U(\mathbb{Z}_{n} \oplus \mathbb{Z}_{m})$, then $\text{deg}(\bar{x})=|U(\mathbb{Z}_{n} \oplus \mathbb{Z}_{m})|$. Thus,\\
			\begin{align*}
				\displaystyle|E|&=\frac{\sum_{\bar{x}\in \mathbb{Z}_{n}\oplus \mathbb{Z}_{m}}\text{deg}(\bar{x})}{2}\\
				& = \frac{\sum_{\bar{x}\in U(\mathbb{Z}_{n} \oplus \mathbb{Z}_{m})}\text{deg}(\bar{x})+\sum_{\bar{x}\in N_{U}(\mathbb{Z}_{n} \oplus \mathbb{Z}_{m})}\text{deg}(\bar{x})}{2}\\
				& = \frac{\phi(m)\phi(n)(\phi(m)\phi(n)-1)
				+[mn-\phi(m)\phi(n)]\phi(m)\phi(n)}{2}\\
				& = \frac{(mn-1)\phi(m)\phi(n)}{2} 	
			\end{align*}
			In other cases,  $2\cdot \bar{x}\notin U(\mathbb{Z}_{n}\oplus\mathbb{Z}_{m}), \ \forall\ \bar{x}\in \mathbb{Z}_{n}\oplus \mathbb{Z}_{m}$. From this, we get $\text{deg}(x)=|U(\mathbb{Z}_{n}\oplus \mathbb{Z}_{m})|, \forall \ \bar{x}\in \mathbb{Z}_{n}\oplus \mathbb{Z}_{m}$. Thus,\\
			\begin{align*}
				\displaystyle|E|&=\frac{\sum_{\bar{x}\in \mathbb{Z}_{n}\oplus \mathbb{Z}_{m}}\text{deg}(\bar{x})}{2}\\
				& = \frac{mn\phi(m)\phi(n)}{2} 	
			\end{align*}
	\end{proof}
\end{theorem}
\begin{lemma}\label{thm24}
	Let $\displaystyle G(\mathbb{Z}_{n}\oplus \mathbb{Z}_{m} )$ be a unit graph. If exactly one of $m$ and $n$ is even, then $\displaystyle G(\mathbb{Z}_{n}\oplus \mathbb{Z}_{m} )$ is bipartite.
	\begin{proof}
		Without loss of generality, assume that $m$ is even and $n$ is odd.\\ Consider the set $W_{1}=\{(x_{1},x_{2})\in \mathbb{Z}_{n}\oplus \mathbb{Z}_{m} \ | \ x_{2}=2\alpha, \ \alpha \in \mathbb{Z} \}$ and\\ $W_{2}=\{(x_{1},x_{2})\in \mathbb{Z}_{n}\oplus \mathbb{Z}_{m} \ | \ x_{2}=2\alpha+1, \ \alpha \in \mathbb{Z} \}$. This sets form separation for $\mathbb{Z}_{n}\oplus \mathbb{Z}_{m}$ such that no vertex in $W_{1}$ is adjacent to any vertex in $W_{1}$. Similarly, for $W_{2}$.
	\end{proof}
\end{lemma}
\begin{corollary}\label{thm25}
	Let $\displaystyle G(\mathbb{Z}_{n} \oplus \mathbb{Z}_{m})$ be a unit graph, where both $m$ and $n$ are even numbers. Then $\displaystyle G(\mathbb{Z}_{n} \oplus \mathbb{Z}_{m})$ is a disconnected graph.
	\begin{proof}
		Let $W_{1}=\{(x_{1},x_{2})\in \mathbb{Z}_{n} \oplus \mathbb{Z}_{m} \ | \ \text{both $x_{1}$ and $x_{1}$ are even (or) both $x_{1}$ and $x_{1}$ are odd } \}$ and $W_{2}=\mathbb{Z}_{n} \oplus \mathbb{Z}_{m} - W_{1}$. If $\bar{x}=(x_{1},x_{2})\in W_{1}$ is adjecent to $\bar{y}=(y_{1},y_{2})\in W_{2}$, then $x_{1}+y_{1}, 
		\ x_{2}+y_{2}$ are units in $\mathbb{Z}_{n}, 
		\mathbb{Z}_{m}$ respectively. But from the construction of $W_{1}$ and $W_{2}$ either $x_{1}+y_{1}$ or $x_{2},y_{2}$ is even, which is contadiction. Hence, no vertex in $W_{1}$ is adjacent to any vertex in $W_{2}$. Hence, $\displaystyle G(\mathbb{Z}_{n} \oplus \mathbb{Z}_{m})$ is a disconnected graph
	\end{proof}
\end{corollary}
\section{Linear Codes from Unit Graph $G(\mathbb{Z}_{p^{n}}\oplus \mathbb{Z}_{q^{m}})$}
In this section, we construct binary and $r$-ary linear codes $C_{2}$ and $C_{r}$ generated from the incidence matrix of the unit graph $\displaystyle G(\mathbb{Z}_{p^{n}}\oplus \mathbb{Z}_{q^{m}})$, where $p$ and $q$ are primes  and $n,\ m\in\mathbb{N}$. We also examine the dual codes $C^{\perp}_{2}$ and $C^{\perp}_{r}$ with their minimum distance.

\begin{theorem}\label{thm21}
Let $\displaystyle G(\mathbb{Z}_{p^{n}} \oplus \mathbb{Z}_{q^{m}})$ be a unit graph, where $\displaystyle p$ and  $\displaystyle q$ are odd primes. Then $G(\mathbb{Z}_{p^{n}} \oplus \mathbb{Z}_{q^{m}})$ is  connected and $\displaystyle \text{diam}(G(\mathbb{Z}_{p^{n}} \oplus \mathbb{Z}_{q^{m}}))\leq 2$.
\begin{proof}
Let $\mathbb{Z}_{p^{n}} \oplus \mathbb{Z}_{q^{m}}=U(\mathbb{Z}_{p^{n}} \oplus \mathbb{Z}_{q^{m}})\cup N_{U}(\mathbb{Z}_{p^{n}} \oplus \mathbb{Z}_{q^{m}})$. We can rewrite
\begin{align}\label{eq0}
	\mathbb{Z}_{p^{n}} \oplus \mathbb{Z}_{q^{m}}&=U(\mathbb{Z}_{p^{n}} \oplus \mathbb{Z}_{q^{m}})\cup N_{(p,u)}\cup N_{(u,q)}\cup N_{(p,q)}\\
\text{Where,} & \nonumber\\
N_{(p,u)} & = \{(x_{1},x_{2})\in \mathbb{Z}_{p^{n}}\oplus\mathbb{Z}_{q^{m}} \ |\ x_{1}=\alpha p \ \& \ x_{2}\in U(\mathbb{Z}_{q^{m}})\}\nonumber\\
N_{(u,q)} & = \{(x_{1},x_{2})\in \mathbb{Z}_{p^{n}}\oplus\mathbb{Z}_{q^{m}} \ | \ x_{1}\in U(\mathbb{Z}_{p^{n}}) \ \& \ x_{2}=\beta q\}\nonumber\\
N_{(p,q)} & = \{(x_{1},x_{2})\in \mathbb{Z}_{p^{n}}\oplus\mathbb{Z}_{q^{m}} \ |\ x_{1}=\alpha p \ \& \ x_{2}=\beta q\}\nonumber
\end{align}
For $\bar{x},\ \bar{y} \in \mathbb{Z}_{p^{n}} \oplus \mathbb{Z}_{q^{m}} $, we have following cases\\
{\bf Case I:} If $\bar{x} , \ \bar{y}\in U(\mathbb{Z}_{p^{n}} \oplus \mathbb{Z}_{q^{m}})$, then $[\bar{x},\bar{0}]$ and $[\bar{0},\bar{y}]$ are edges in $\displaystyle G(\mathbb{Z}_{p^{n}} \oplus \mathbb{Z}_{q^{m}})$.\\
{\bf Case II:}  If $\bar{x} \in U(\mathbb{Z}_{p^{n}} \oplus \mathbb{Z}_{q^{m}})$ and $\bar{y}\in N(\mathbb{Z}_{p^{n}} \oplus \mathbb{Z}_{q^{m}})$, then $\bar{x}=(x_{1},x_{2})$ and for $\bar{y}$ we have following subcases:\\
(a) If $\bar{y}\in N_{(p,u)}$, then  $\bar{y}=(\alpha p,u_{1})$. If $\bar{x}$ and $\bar{y}$ are adjacent, then $d(x,y)=1\leq 2$. Suppose, $\bar{x}$ and $\bar{y}$ are not adjacent. Then, there exist $\bar{z}=(x_{1},0)\in N_{(u,q)}$ such that $[\bar{x},\bar{z}]$ and   $[\bar{z},\bar{x}]$ are edges in $\displaystyle G(\mathbb{Z}_{p^{n}} \oplus \mathbb{Z}_{q^{m}})$, this gives $d(\bar{x},\bar{y})\leq 2$.\\
(b) If $\bar{y}\in N_{(p,q)}$ then $\bar{y}=(\alpha p,\beta q)$. Clearly, $\bar{x}$ and $\bar{y}$ are adjacent.\\
{\bf Case III:}  If $\bar{x},\bar{y}\in N(\mathbb{Z}_{p^{n}} \oplus \mathbb{Z}_{q^{m}})$, then following possibilities arise:\\
(a) If $\bar{x}, \bar{y} \in N_{(p,u)}$, then $\bar{x}=(\alpha_{1}p,u_{1})$ and $\bar{y}=(\alpha_{2}p,u_{2})$.
Consider $\bar{z}=(u_{3},0)\in N_{(u,q)}$, from this, we get, $[\bar{x},\bar{z}]$ and $[\bar{z},\bar{y}]$ are edges in $\displaystyle G(\mathbb{Z}_{p^{n}} \oplus \mathbb{Z}_{q^{m}})$. Hence $d(\bar{x},\bar{y})\leq 2$. Similarly, if $\bar{x}$ and $\bar{y}$ are in same set in (\ref{eq0}), then there is $\bar{z}$ in any other set in (\ref{eq0}) such that, $\bar{x}$ and $\bar{y}$ are adjacent to $\bar{z}$.\\
(b) If $\bar{x} \in N_{(p,u)}$ and  $\bar{y} \in N_{(u,q)}$, then
	$\bar{x}=(\alpha p,u_{1})$ and $\bar{y}=(u_{2},\beta q)$. Clearly, $\bar{x}$ and $\bar{y}$ are adjacent, which implies, $d(\bar{x},\bar{y})=1\leq 2$.\\
(c) If $\bar{x} \in N_{(p,q)}$ and  $\bar{y} \in N_{(u,q)}$, then $\bar{x}=(\alpha p, \beta_{1}q)$ and $\bar{y}=(u_{1}, \beta_{2}q)$. Consider, the element, $\bar{z}=(u_{1},u_{2})\in U(\mathbb{Z}_{p^{n}} \oplus \mathbb{Z}_{q^{m}})$, we get  $\bar{x}$ is adjacent to $\bar{z}$ and $\bar{z}$ is  adjacent to $\bar{y}$. Hence, $d(\bar{x},\bar{y})\leq 2$. \\
All other cases follows in the same manner.\\
Hence, $d(\bar{x}, \bar{y})\leq 2$, for all $\bar{x},\bar{y}\in \mathbb{Z}_{p^{n}} \oplus \mathbb{Z}_{q^{m}}$ and $G(\mathbb{Z}_{p^{n}} \oplus \mathbb{Z}_{q^{m}})$ connected graph, which
gives diam$(G(\mathbb{Z}_{p^{n}} \oplus \mathbb{Z}_{q^{m}}))\leq 2$.
\end{proof}
\end{theorem}
\begin{corollary}\label{thm23}
	Let $\displaystyle G(\mathbb{Z}_{p^{n}} \oplus \mathbb{Z}_{q^{m}})$ be a unit graph, where $\displaystyle p$ and  $\displaystyle q$ are odd primes. Then $\displaystyle \lambda(G(\mathbb{Z}_{p^{n}} \oplus \mathbb{Z}_{q^{m}}))=\phi(p^{n})\phi(q^{m})-1$.
\begin{proof}
It follows from Theorem \ref{thm3} and \ref{thm21}.
\end{proof}
\end{corollary}

\begin{theorem}\label{thm26}
	 Let $\displaystyle G(\mathbb{Z}_{p^{n}} \oplus \mathbb{Z}_{2^{m}})$ be a unit graph, where $\displaystyle p$ be an odd prime. Then $\displaystyle G(\mathbb{Z}_{p^{n}} \oplus \mathbb{Z}_{2^{m}})$ is a connected bipartite graph and $\displaystyle \lambda(G(\mathbb{Z}_{p^{n}} \oplus \mathbb{Z}_{2^{m}}))=2^{m-1}\phi(p^{n})$.
\begin{proof}
We know that $\mathbb{Z}_{2^{m}p^{n}}\simeq\mathbb{Z}_{p^{n}}\oplus \mathbb{Z}_{2^{m}}$. From this we have, $G(\mathbb{Z}_{2^{m}p^{n}})$ is isomorphic to $G(\mathbb{Z}_{p^{n}}\oplus \mathbb{Z}_{2^{m}})$. Hence result follows from, Theorem \ref{thm9} and Lemma \ref{thm24}.
\end{proof}
\end{theorem}

\begin{theorem}
Let $\displaystyle G(\mathbb{Z}_{p^{n}} \oplus \mathbb{Z}_{q^{m}})$ be a unit graph and $H$ be a $|V|\times |E|$ incidence matrix of $\displaystyle G(\mathbb{Z}_{p^{n}} \oplus \mathbb{Z}_{q^{m}})$.
\begin{enumerate}
  \item If both $p$ and $q$ are odd primes then\\ $\displaystyle C_{2}(H)=\left[\frac{(p^{n}q^{m}-1)\phi(p^{n})\phi(q^{m})}{2},p^{n}q^{m}-1,\phi(p^{n})\phi(q^{m})-1\right]_{2}$ is the binary code generated by $H$ over the finite field $\displaystyle \mathbb{F}_{2}$.
  \item If $p$ is odd prime and $q$ is even prime, then for any odd prime $r$, $\displaystyle C_{r}(H)=[p^{n}\phi(p^{n})2^{2(m-1)},2^{m}p^{n}-1,2^{m-1}\phi(p^{n})]_{r}$ is the $r$-ary code generated by $H$ over the finite field $\displaystyle \mathbb{F}_{r}$.
\end{enumerate}
\begin{proof}
\begin{enumerate}
  \item Let $\displaystyle G(\mathbb{Z}_{p^{n}} \oplus \mathbb{Z}_{q^{m}})$ be a unit graph, where $p$ and $q$ both are odd primes and $H$ be an incidence matrix of $\displaystyle G(\mathbb{Z}_{p^{n}} \oplus \mathbb{Z}_{q^{m}})$.  By Theorem \ref{thm21},  $\displaystyle G(\mathbb{Z}_{p^{n}} \oplus \mathbb{Z}_{q^{m}})$ is connected graph and hence by Theorem \ref{thm2}, binary code generated by $H$ is $C_{2}(H)=[|E|,|V|-1,\lambda(\displaystyle G(\mathbb{Z}_{p^{n}} \oplus \mathbb{Z}_{q^{m}}))]_{2}$. By, Theorem \ref{thm41} and Corollary \ref{thm23},\\ we get
$\displaystyle |E|=\frac{\phi(p^{n})\phi(q^{m})(p^{n}q^{m}-1)}{2}, \ |V|=p^{n}q^{m}$  and the  edge connectivity of $G(\mathbb{Z}_{p^{n}} \oplus \mathbb{Z}_{q^{m}})$ is $\lambda(G(\mathbb{Z}_{p^{n}} \oplus \mathbb{Z}_{q^{m}}))=\phi(p^{n})\phi(q^{m})-1$.\\ Hence we get $\displaystyle C_{2}(H)=\left[\frac{(p^{n}q^{m}-1)\phi(p^{n})\phi(q^{m})}{2},p^{n}q^{m}-1,\phi(p^{n})\phi(q^{m})-1\right]_{2}$
\item Let $q=2$ and $p$ be an odd prime. Then  by Theorem \ref{thm26}, $\displaystyle G(\mathbb{Z}_{p^{n}} \oplus \mathbb{Z}_{2^{m}})$  is a connected bipartite graph and hence by Theorem \ref{thm5}, for any odd prime $r$, $r$-ary code generated by $H$ is $C_{r}(H)=[|E|,|V|-1,\lambda(\displaystyle G(\mathbb{Z}_{p^{n}} \oplus \mathbb{Z}_{2^{m}}))]_{q}$. Using, Theorem \ref{thm26}, we get
    $\displaystyle |E|=p^{n}\phi(p^{n})2^{2(m-1)}, \ |V|=p^{n}2^{m}$  and the  edge connectivity of $\displaystyle G(\mathbb{Z}_{p^{n}} \oplus \mathbb{Z}_{2^{m}})$ is $\lambda(\displaystyle G(\mathbb{Z}_{p^{n}} \oplus \mathbb{Z}_{2^{m}})=2^{m-1}\phi(p^{n})$. Hence,
    $\displaystyle C_{r}(H)=[2^{2(m-1)}p^{n}\phi(p^{n}),2^{m}p^{n}-1,2^{m-1}\phi(p^{n})]_{r}$.
\end{enumerate}
\end{proof}
\end{theorem}

\begin{corollary}
Let $C_{r}(H)$ and $C_{2}(H)$ denote the codes generated by incidence matrix of $\displaystyle G(\mathbb{Z}_{p^{n}}\oplus\mathbb{Z}_{2^{m}})$ and $\displaystyle G(\mathbb{Z}_{p^{n}}\oplus \mathbb{Z}_{q^{m}})$. Then
\begin{enumerate}
  \item Dual of code $C_{2}(H)$ is $\displaystyle C^{\perp}_{2}=\left[ \frac{(p^{n}q^{m}-1)\phi(p^{n})\phi(q^{m})}{2}, \frac{(p^{n}q^{m}-1)[\phi(p^{n})\phi(q^{m})-2]}{2}, 3
   \right]_{2}$.
   \item Dual of code $C_{r}(H)$ is $\displaystyle C^{\perp}_{r}=\left[p^{n}\phi(p^{n})2^{2(m-1)},p^{n}2^{m}(\phi(p^{n})2^{m-2}-1)+1,4\right]_{q}$, for $2^{m}p^{n}\neq 6$.
\end{enumerate}
\begin{proof}
\begin{enumerate}
  \item From Theorems \ref{thm8}, \ref{thm6}, we have dimension of $C_{2}^{\perp}$ is $\frac{(p^{n}q^{m}-1)[\phi(p^{n})\phi(q^{m})-2]}{2}$ and minimum distance is $g_{r}(G(\mathbb{Z}_{p^{n}} \oplus \mathbb{Z}_{q^{m}}))$.
  Note that $\bar{x}=(u_{1},\beta q)\in N_{(u,q)}$ and $\bar{y}=(\alpha p,u_{2})\in N_{(p,u)}$ are adjacent and both $\bar{x}, \ \bar{y}$ are adjacent to $\bar{z}=(u_{1},u_{2})\in U(\mathbb{Z}_{p^{n}} \oplus \mathbb{Z}_{q^{m}})$. Hence, we get $\displaystyle g_{r}(G(\mathbb{Z}_{p^{n}} \oplus \mathbb{Z}_{q^{m}}))=3$.
  \item This result follows from Theorems \ref{thm8}, \ref{thm6} and \ref{thm10}.
\end{enumerate}
\end{proof}
\end{corollary}

\section{Linear Codes from the Unit Graph $\displaystyle G(\mathbb{Z}_{p_{1}^{n_1}p^{n_2}_{2}} \oplus \mathbb{Z}_{q_{1}^{m_1}q^{m_2}_{2}})$}
In this section, we extend the results in Section 4 and construct binary and $r$-ary linear codes from the incidence matrix of the unit graph $\displaystyle G(\mathbb{Z}_{p_{1}^{n_1}p^{n_2}_{2}} \oplus \mathbb{Z}_{q_{1}^{m_1}q^{m_2}_{2}})$, where $p_{1},p_{2},q_{1}$ and $q_{2}$  are primes and $\displaystyle n_{1},n_{2},m_{1},m_{2}\in \mathbb{N}$.

\begin{theorem}\label{thm31}
Let $G(\mathbb{Z}_{p_{1}^{n_1}p^{n_2}_{2}} \oplus \mathbb{Z}_{q_{1}^{m_1}q^{m_2}_{2}})$ be a unit graph, where $\displaystyle p_{1}, p_{2},q_{1}$ and $q_{2}$ are odd primes. Then $G(\mathbb{Z}_{p_{1}^{n_1}p^{n_2}_{2}} \oplus \mathbb{Z}_{q_{1}^{m_1}q^{m_2}_{2}})$ is a connected graph and $\displaystyle \text{diam}(G(\mathbb{Z}_{p_{1}^{n_1}p^{n_2}_{2}} \oplus \mathbb{Z}_{q_{1}^{m_1}q^{m_2}_{2}}))\leq 2$.
\begin{proof}
 We can rewrite
 \begin{align}\label{eq1}
 \displaystyle\mathbb{Z}_{p_{1}^{n_1}p^{n_2}_{2}} \oplus \mathbb{Z}_{q_{1}^{m_1}q^{m_2}_{2}}&=U(\mathbb{Z}_{p_{1}^{n_1}p^{n_2}_{2}} \oplus \mathbb{Z}_{q_{1}^{m_1}q^{m_2}_{2}})\cup N_{U}(\mathbb{Z}_{p_{1}^{n_1}p^{n_2}_{2}} \oplus \mathbb{Z}_{q_{1}^{m_1}q^{m_2}_{2}})\nonumber\\
 & = U(\mathbb{Z}_{p_{1}^{n_1}p^{n_2}_{2}} \oplus \mathbb{Z}_{q_{1}^{m_1}q^{m_2}_{2}})\bigcup_{i=1}^{2} N_{(p_{i},u)}\bigcup_{i=1}^{2} N_{(u,q_{i})}\bigcup_{i=1}^{2}\bigcup_{j=1}^{2} N_{(p_{i},q_{j})}\bigcup_{i=1}^{2} N_{(p_{i},q_{1}q_{2})}\notag\\
 & \ \ \bigcup_{i=1}^{2} N_{(p_{1}p_{2},q_{i})}\bigcup N_{(p_{1}p_{2},q_{1}q_{2})}\bigcup N_{(p_{1}p_{2},u)}\bigcup N_{(u,q_{1}q_{2})}
 \end{align}
 where
 \begin{align*}
 N_{(p_{i},u)}&=\{(x_{1},x_{2})\in \mathbb{Z}_{p_{1}^{n_1}p^{n_2}_{2}} \oplus \mathbb{Z}_{q_{1}^{m_1}q^{m_2}_{2}} \ | \ x_{2}\in U(\mathbb{Z}_{q_{1}^{m_1}q^{m_2}_{2}}), \ x_{1}=\alpha p_{i} \ \& \ p_{j}\nmid \alpha, \text{ for } \ i\neq j \}\\
N_{(u,q_{i})}&=\{(x_{1},x_{2})\in \mathbb{Z}_{p_{1}^{n_1}p^{n_2}_{2}} \oplus \mathbb{Z}_{q_{1}^{m_1}q^{m_2}_{2}} \ | \ x_{1}\in U(\mathbb{Z}_{p_{1}^{m_1}p^{m_2}_{2}}),\ x_{2}=\beta q_{i} \ \& \ q_{j}\nmid \beta, \text{ for } i\neq j  \}\\
 N_{(p_{i},q_{j})}&=\{(x_{1},x_{2})\in \mathbb{Z}_{p_{1}^{n_1}p^{n_2}_{2}} \oplus \mathbb{Z}_{q_{1}^{m_1}q^{m_2}_{2}} \ | \ x_{1}=\alpha p_{i} \ \& \ p_{k}\nmid \alpha \text{ for } k\neq i, \notag\\ & \quad x_{2}=\beta q_{j} \ \& \ q_{l}\nmid \beta, \text{ for } j\neq l  \}\\
  N_{(p_{i},q_{1}q_{2})}&=\{(x_{1},x_{2})\in \mathbb{Z}_{p_{1}^{n_1}p^{n_2}_{2}} \oplus \mathbb{Z}_{q_{1}^{m_1}q^{m_2}_{2}} \ | \ x_{1}=\alpha p_{i} \ \& \ p_{j}\nmid \alpha \text{ for } j\neq i, \  x_{2}=\beta q_{1}q_{2}  \}\\
   N_{(p_{1}p_{2},q_{i})}&=\{(x_{1},x_{2})\in \mathbb{Z}_{p_{1}^{n_1}p^{n_2}_{2}} \oplus \mathbb{Z}_{q_{1}^{m_1}q^{m_2}_{2}} \ | \ x_{1}=\alpha p_{1}p_{2},\  x_{2}=\beta q_{i} \ \& \ q_{j}\nmid \beta, \text{ for } i\neq j  \}\\
    N_{(p_{1}p_{2},q_{1}q_{2})}&=\{(x_{1},x_{2})\in \mathbb{Z}_{p_{1}^{n_1}p_{2}^{n_2}} \oplus \mathbb{Z}_{q_{1}^{m_1}q_{2}^{m_2}} \ | \ x_{1}=\alpha p_{1}p_{2}, \ x_{2}=\beta q_{1}q_{2}   \}\\
 N_{(p_{1}p_{2},u)}&=\{(x_{1},x_{2})\in \mathbb{Z}_{p_{1}^{n_1}p_{2}^{n_2}} \oplus \mathbb{Z}_{q_{1}^{m_1}q_{2}^{m_2}} \ | \ x_{1}=\alpha p_{1}p_{2}, \ x_{2}\in U(\mathbb{Z}_{q_{1}^{m_1}q_{2}^{m_2}})  \}\\
 N_{(u,q_{1}q_{2})}&=\{(x_{1},x_{2})\in \mathbb{Z}_{p_{1}^{n_1}p_{2}^{n_2}} \oplus \mathbb{Z}_{q_{1}^{m_1}q_{2}^{m_2}} \ |  \ x_{1}\in U(\mathbb{Z}_{p_{1}^{n_1}p_{2}^{n_2}}), \ x_{2}=\beta q_{1}q_{2},  \}
     \end{align*}
 For $\displaystyle \bar{x},\bar{y}\in \mathbb{Z}_{p_{1}^{n_1}p^{n_2}_{2}} \oplus \mathbb{Z}_{q_{1}^{m_1}q^{m_2}_{2}}$, consider  following cases:\\
{\bf Case I:} If $\displaystyle \bar{x},\bar{y}\in U(\mathbb{Z}_{p_{1}^{n_1}p^{n_2}_{2}} \oplus \mathbb{Z}_{q_{1}^{m_1}q^{m_2}_{2}})$, then $\bar{x}$ is adjacent to $\bar{0}$ and $\bar{0}$ is adjacent to $\bar{y}$, this gives $d(\bar{x},\bar{y})\leq 2$\\
{\bf Case II:} If $\displaystyle \bar{x}\in U(\mathbb{Z}_{p_{1}^{n_1}p^{n_2}_{2}} \oplus \mathbb{Z}_{q_{1}^{m_1}q^{m_2}_{2}})$ and $\displaystyle \bar{y}\in N_{U}(\mathbb{Z}_{p_{1}^{n_1}p^{n_2}_{2}} \oplus \mathbb{Z}_{q_{1}^{m_1}q^{m_2}_{2}})$, then $\bar{x}=(x_{1},x_{2})$ and for $\bar{y}$ we have following possibilities:\\
(a) If $\displaystyle \bar{y}\in N_{(p_{i},u)}$ ,then $y=(\alpha p_{i},u)$. Suppose $[\bar{x},\bar{y}]$ is not edge in $G(\mathbb{Z}_{p_{1}^{n_1}p^{n_2}_{2}} \oplus \mathbb{Z}_{q_{1}^{m_1}q^{m_2}_{2}})$. Then, either $x_{1}+\alpha p_{i}\notin U(\mathbb{Z}_{p^{n_{1}}_{1}}\oplus \mathbb{Z}_{p^{n_{2}}_{2}})$ or $x_{2}+u\notin U(\mathbb{Z}_{q^{m_{1}}_{1}\oplus \mathbb{Z}_{q^{m_{2}}_{2}}})$.\\
If $x_{1}+\alpha p_{i}\notin U(\mathbb{Z}_{p^{n_{1}}_{1}}\oplus \mathbb{Z}_{p^{n_{2}}_{2}})$ and $x_{2}+u\in U(\mathbb{Z}_{q^{m_{1}}_{1}\oplus \mathbb{Z}_{q^{m_{2}}_{2}}})$, then $p_{j}\mid x_{1}+\alpha p_{i}$, which gives, $x_{1}+\alpha p_{i}=\alpha_{2}p_{j}$, for $i\neq j$. Consider, $\bar{z}=(\alpha_{2}p_{j},u)\in N_{(p_{j},u)}$. Then $[\bar{x},\bar{z}]$ and $[\bar{z},\bar{y}]$ are edges in $G(\mathbb{Z}_{q^{m_{1}}_{1}}\oplus \mathbb{Z}_{q^{m_{2}}_{2}})$.\\
If $x_{1}+\alpha p_{i}\in U(\mathbb{Z}_{p^{n_{1}}_{1}}\oplus \mathbb{Z}_{p^{n_{2}}_{2}})$ and $x_{2}+u\notin U(\mathbb{Z}_{q^{m_{1}}_{1}\oplus \mathbb{Z}_{q^{m_{2}}_{2}}})$, then $\bar{x}$ is adjacent to $\bar{z}=(x_{1},0)$ and $\bar{z}=(x_{1},0)$ is adjacent to $\bar{y}$.\\
If $x_{1}+\alpha p_{i}\notin U(\mathbb{Z}_{p^{n_{1}}_{1}}\oplus \mathbb{Z}_{p^{n_{2}}_{2}})$ and $x_{2}+u\notin U(\mathbb{Z}_{q^{m_{1}}_{1}\oplus \mathbb{Z}_{q^{m_{2}}_{2}}})$, then $\bar{z}=(\alpha_{2}p_{j},0)\in N_{(p_{j},q_{1}q_{2})}$, where $\alpha_{2}p_{j}=x+\alpha_{1}p_{i}$ is adjacent to both $\bar{x}$ and $\bar{y}$. Hence, $d(\bar{x},\bar{y})\leq 2$.\\
(b) If $\displaystyle \bar{y}\in N_{(p_{i},q_{j})}$, then $\bar{y}=(\alpha_{1}p_{i},\beta_{1}q_{j})$. If $\bar{x}$ is not adjacent to
$\bar{y}$, then either $p_{k}\mid x_{1}+\alpha_{1}p_{i}$ or $q_{l}\mid x_{2}+\beta_{1}q_{j}$, for $i\neq k$ and $j\neq l$.\\
If $p_{k}\mid x_{1}+\alpha_{1}p_{i}$ and $q_{l}\nmid x_{2}+\beta_{1}q_{j}$. Consider, $\bar{z}=(\alpha_{2}p_{k},x_{2})\in N_{(p_{k},u)}$, where $\alpha_{2}p_{k}=x_{1}+\alpha_{1}p_{i}$. Clearly, $[\bar{x},\bar{z}]$ and $[\bar{z},\bar{y}]$ are edges in $G(\mathbb{Z}_{p_{1}^{n_{1}}p_{2}^{n_{2}}}\oplus\mathbb{Z}_{q_{1}^{m_{1}}q_{2}^{m_{2}}})$. Other two possibilities follows in the same way. Hence, we get $d(\bar{x},\bar{y})\leq 2$.\\
(c)  If $\displaystyle \bar{y}\in N_{(p_{i},q_{1}q_{2})}$, then $\bar{y}=(\alpha_{1}p_{i},\beta_{1}q_{1}q_{2})$. If $\bar{x}$ and $\bar{y}$ are not adjacent, then $p_{j}\mid x_{1}+\alpha_{1}p_{i}$, for $i\neq j$. Now, consider $\bar{z}=(\alpha_{2}p_{j},x_{2})\in N_{(p_{j},u)}$, where $\alpha_{2}p_{j}=x_{1}+\alpha_{1}p_{i}$. Then $\bar{x}$ is adjacent to $
\bar{z}$ and $\bar{z}$ is adjacent to $\bar{y}$, which gives $d(\bar{x},\bar{y})\leq 2$.\\
(d) If $\displaystyle \bar{y}\in N_{(p_{1}p_{2},u)}$, then $\bar{y}=(\alpha_{1}p_{1}p_{2},u_{1})$. Suppose, $[\bar{x},\bar{y}]$ is not an edge in $G(\mathbb{Z}_{p_{1}^{n_1}p_{2}^{n_2}}\oplus \mathbb{Z}_{q_{1}^{m_1}q_{2}^{m_2}} )$.  Take an element $\bar{z}=(x_{1},0)\in N_{(u,q_{1}q_{2})}$. Then $[\bar{x},\bar{z}]$ and $[\bar{z},\bar{y}]$ are edges in $G(\mathbb{Z}_{p_{1}^{n_1}p_{2}^{n_2}}\oplus \mathbb{Z}_{q_{1}^{m_1}q_{2}^{m_2}})$, which implies $d(\bar{x},\bar{y})\leq 2$.\\
Using the above procedures  for all other possibilities in this case, we can prove  that $d(\bar{x},\bar{y})\leq 2$.\\ 
{\bf Case III:} Let $\displaystyle \bar{x}, \bar{y}\in N_{U}(\mathbb{Z}_{p_{1}^{n_1}p^{n_2}_{2}} \oplus \mathbb{Z}_{q_{1}^{m_1}q^{m_2}_{2}})$. We prove for some possibilities for $\bar{x}$ and $\bar{y}$, we give proof for some possibilities and all other possibilities follows from the same procedure:\\
(a) If $\bar{x},\bar{y}\in N_{(p_{i},u)}$, then for any $\bar{z}=(\alpha_{3}p_{j},0)\in N_{(p_{j},q_{1}q_{2})}$, for $i\neq j$, we have the edges in $[\bar{x},\bar{z}]$ and $[\bar{z},\bar{y}]$ in $G(\mathbb{Z}_{p_{1}^{n_1}p_{2}^{n_2}}\oplus\mathbb{Z}_{q_{1}^{m_1}q_{2}^{m_2}})$. Similarly, to this if $\bar{x},\ \bar{y}$ are in same sets in (\ref{eq1}), then there exist $\bar{z}$ in any other set in (\ref{eq1}) such that both $\bar{x}, \bar{y}$ are adjacent to $\bar{z}$.\\
(b) If $\bar{x}\in N_{(p_{i},u)}$ and $\bar{y}\in N_{(u,q_{j})}$, then $\bar{x}=(\alpha_{1}p_{i},u_{1})$ and $\bar{y}=(u_{2},\beta_{1}q_{j})$. Suppose, if $\bar{x}$ is not adjacent to $\bar{y}$, then either $p_{k}\mid\alpha_{1}p_{i}+u_{2}$, for $k\neq i$ or $q_{l}\mid u_{1}+\beta_{1}q_{j}$, for $j\neq l$. If $p_{k}\mid\alpha_{1}p_{i}+u_{2}$ and  $q_{l}\nmid u_{1}+\beta_{1}q_{j}$, for $j\neq l$, then take $\bar{z}=(\alpha_{2}p_{k},u_{1})\in N_{(p_{k},u)}$, where $\alpha_{2}p_{k}=\alpha_{1}p_{i}+u_{2}$. Clearly, $\bar{x}$ and $\bar{y}$ are adjacent to $\bar{z}$. Similarly, for other possibility, we have,  $d(\bar{x},\bar{y})\leq 2$.\\
(c) If $\bar{x}\in N_{(p_{i},u)}$ and $\bar{y}\in N_{(p_{j},q_{k})}$, then $\bar{x}=(\alpha_{1}p_{i},u_{1})$ and $\bar{y}=(\alpha_{1}p_{j},\beta_{1}q_{k})$. If $i\neq j$ and $q_{l}\mid u_{1}+\beta_{1}q_{k}$, then for $\bar{z}=(\alpha_{1}p_{i}+\alpha_{2}p_{j},\beta_{2}q_{l})$, where, $q_{l}\beta_{2}=u_{1}+\beta_{1}q_{k}$, we get $[\bar{x},\bar{z}]$ and $[\bar{z},\bar{y}]$ are edges in $G(\mathbb{Z}_{p_{1}^{n_{1}}p_{2}^{n_2}}\oplus\mathbb{Z}_{q_{1}^{m_2}q_2^{m_2}})$. For $i=j$ and if $q_{l}\mid \beta_{1}q_{k}+u_{1}$, then take $\bar{z}=(\alpha_{3}p_{v},\beta_{2}q_{l})$, where $v\neq i,j$ and $\beta_{2}q_{l}=u_{1}+\beta_{1}q_{k}$. Hence, we get $d(\bar{x},\bar{y})\leq 2$.\\
(d) If $\bar{x}\in N_{(p_{i},u)}$ and $\bar{y}\in N_{(p_{1}p_{2},q_{j})}$, then $\bar{x}=(\alpha_{1}p_{i},u_{1})$ and $\bar{y}=(\alpha_{2}p_{1}p_{2},\beta_{1}q_{j})$. If $q_{l}\mid u_{1}+\beta_{1}q_{j}$, then take an element $\bar{z}=(\alpha_{1}p_{i}+\alpha_{3}p_{k},\beta_{2}q_{l})$ is adjacent to both $\bar{x}$ and $\bar{y}$, where $p_{i}\nmid \alpha_{3} \ \& \ \beta_{2}q_{l}=u_{1}+\beta_{1}q_{j} $.\\
(e) If $\bar{x}\in N_{(p_{i},u)}$ and $\bar{y}\in N_{(p_{1}p_{2},u)}$, then $\bar{x}=(\alpha_{1}p_{i},u_{1})$ and $\bar{y}=(\alpha_{2}p_{1}p_{2},u_{2})$. Take $\bar{z}=(\alpha_{1}p_{i}+\alpha_{2}p_{j},0)$, where $i\neq j$ and $p_{i}\nmid \alpha_{2}$. From this, we get $[\bar{x},\bar{z}]$ and $[\bar{z},\bar{y}]$ are edges in $G(\mathbb{Z}_{p_{1}^{n_1}p_{2}^{n_2}}\oplus\mathbb{Z}_{q_{1}^{m_1}q_{2}^{m_2}})$. Hence, $d(\bar{x},\bar{y})\leq 2$.\\
(f) If $\bar{x}\in N_{(p_{i},q_{j})}$ and $\bar{y}\in N_{(p_{1}p_{2},q_{k})}$, then $\bar{x}=(\alpha_{1}p_{i},\beta_{1}q_{j})$ and $\bar{y}=(\alpha_{2}p_{1}p_{2},\beta_{2}q_{k})$. If $j=k$, then take $\bar{z}=(\alpha_{1}p_{i}+\alpha_{3}p_{l},\beta_{3}q_{v})\in N_{(u,q_{v})}$, where $l\neq i, \ v\neq j,k $. Clearly $\bar{x}$ and $\bar{y}$ are adjacent to $\bar{z}$. If $j\neq k$, then take $\bar{z}=(\alpha_{1}p_{i}+\alpha_{3}p_{l},\beta_{1}q_{j}+\beta_{2}q_{k})\in U(\mathbb{Z}_{p_{1}^{n_1}p_{2}^{n_2}}\oplus \mathbb{Z}_{q_{1}^{m_1}q_{2}^{m_2}})$, where $l\neq i \ \& \ j\neq k $, from this we get $[\bar{x},\bar{z}]$ and $[\bar{z},\bar{y}]$ are edges in $G(\mathbb{Z}_{p_{1}^{n_1}p_{2}^{n_2}}\oplus\mathbb{Z}_{q_{1}^{m_1}q_{2}^{m_2}})$. This gives $d(\bar{x},\bar{y})\leq 2$\\
(g) If $\bar{x}\in N_{(p_{i},q_{j})}$ and $\bar{y}\in N_{(p_{1}p_{2},u)}$, then $\bar{x}=(\alpha_{1}p_{i},\beta_{1}q_{j})$ and $\bar{y}=(\alpha_{2}p_{1}p_{2},u_{1})$. If $q_{k}\mid u_{1}+\beta_{1}q_{j} $, then take $\bar{z}=(\alpha_{1}p_{i}+\alpha_{3}p_{k},\beta_{2}q_{k})$, where $i\neq k, p_{i}\nmid \alpha_{3}$ and $\beta_{2}q_{k}=u_{1}+\beta_{1}q_{j}$. If $q_{k}\nmid u_{1}+\beta_{1}q_{j} $, then take $\bar{z}=(\alpha_{1}p_{i},\alpha_{2}p_{j},u_{1})$. Hence, $d(\bar{x},\bar{y})\leq 2$.\\
(h) If $\bar{x}\in N_{(p_{i},q_{1}q_{2})}$ and $\bar{y}\in N_{(p_{1}p_{2},q_{j})}$, then $\bar{x}=(\alpha_{1}p_{i},\beta_{1}q_{1}q_{2})$ and $\bar{y}=(\alpha_{2}p_{1}p_{2},\beta_{2}q_{j})$. Take $\bar{z}=(\alpha_{1}p_{i}+\alpha_{3}p_{k},\beta_{2}q_{j}+\beta_{3}q_{l})\in U(\mathbb{Z}_{p_{1}^{n_1}p_{2}^{n_2}}\oplus\mathbb{Z}_{q_{1}^{m_1}q_{2}^{m_2}})$, for $i\neq k$ and $j\neq l$. From this, we get $\bar{x}$ and $\bar{y}$ are adjacent to $\bar{z}$ and hence $d(\bar{x},\bar{y})\leq 2$.\\
(i) If $\bar{x}\in N_{(p_{i},q_{1}q_{2})}$ and $\bar{y}\in N_{(p_{1}p_{2},q_{1}q_{2})}$, then $\bar{x}=(\alpha_{1}p_{i},\beta_{1}q_{1}q_{2})$ and $\bar{y}=(\alpha_{2}p_{1}p_{2},\beta_{2}q_{1}q_{2})$. Take $\bar{z}=(\alpha_{1}p_{i}+\alpha_{3}p_{k},u)$, for $i\neq k, \ p_{i}\nmid \alpha_{3}$ and $u\in U(\mathbb{Z}_{q_{1}^{m_1}q_{2}^{m_2}})$. Then, $\bar{x}$ and $\bar{y}$ are adjacent to $\bar{z}$ and hence $d(\bar{x},\bar{y})\leq 2$.\\
 Hence,  $d(\bar{x}, \bar{y})\leq 2$, for all $\bar{x},\bar{y}\in \mathbb{Z}_{p_{1}^{n_1}p_{2}^{n_2}} \oplus \mathbb{Z}_{q_{1}^{m_1}q_{2}^{m_2}}$ and $G(\mathbb{Z}_{p_{1}^{n_1}p_{2}^{n_2}} \oplus \mathbb{Z}_{q_{1}^{m_1}q_{2}^{m_2}})$ is a connected graph, which
gives diam$(G(\mathbb{Z}_{p_{1}^{n_1}p_{2}^{n_2}} \oplus \mathbb{Z}_{q_{1}^{m_1}q_{2}^{m_2}}))\leq 2$.
\end{proof}
\end{theorem}
\begin{corollary}\label{thm33}
	Let $\displaystyle G(\mathbb{Z}_{p_{1}^{n_1}p_{2}^{n_2}} \oplus \mathbb{Z}_{q_{1}^{m_1}q_{2}^{m_2}})$ be a unit graph, where $\displaystyle p_{1},p_{2},q_{1}$ and  $\displaystyle q_{2}$ are odd primes. Then $\displaystyle \lambda(G(\mathbb{Z}_{p_{1}^{n_1}p_{2}^{n_2}} \oplus \mathbb{Z}_{q_{1}^{m_1}q_{2}^{m_2}}))=\phi(p_{1}^{n_1}p_{2}^{n_2})\phi(q_{1}^{m_1}q_{2}^{m_2})-1$.
\begin{proof}
This result follows from Theorem \ref{thm3} and \ref{thm31}
\end{proof}
\end{corollary}

\begin{theorem}\label{thm34}
	Let $\displaystyle G(\mathbb{Z}_{p_{1}^{n_1}p_{2}^{n_2}}\oplus\mathbb{Z}_{2^{m_1}q^{m_{2}}})$ be a unit graph, where $p_{1},p_{2}$ and $q_{1}$ are odd primes. Then $\displaystyle G(\mathbb{Z}_{p_{1}^{n_1}p_{2}^{n_2}}\oplus\mathbb{Z}_{2^{m_1}q^{m_{2}}})$ is a connected graph and $\displaystyle \text{diam}( G(\mathbb{Z}_{p_{1}^{n_1}p_{2}^{n_2}}\oplus\mathbb{Z}_{2^{m_1}q^{m_{2}}}))\leq 3$.
	\begin{proof}
		As in the Theorem \ref{thm31}, we rewrite $\mathbb{Z}_{p_{1}^{n_1}p_{2}^{n_2}}\oplus\mathbb{Z}_{2^{m_1}q^{m_{2}}}$ as follows,
		\begin{align}\label{eq2}
		 \displaystyle \mathbb{Z}_{p_{1}^{n_1}p_{2}^{n_2}}\oplus\mathbb{Z}_{2^{m_1}q^{m_{2}}}&=U(\mathbb{Z}_{p_{1}^{n_1}p_{2}^{n_2}}\oplus\mathbb{Z}_{2^{m_1}q^{m_{2}}})\cup N_{U}(\mathbb{Z}_{p_{1}^{n_1}p_{2}^{n_2}}\oplus\mathbb{Z}_{2^{m_1}q^{m_{2}}})\nonumber\\ &=U(\mathbb{Z}_{p_{1}^{n_1}p_{2}^{n_2}}\oplus\mathbb{Z}_{2^{m_1}q^{m_{2}}})\bigcup_{i=1}^{2}N_{(p_{i},u)}\bigcup_{i=1}^{2}N_{(p_{i},2)}\bigcup_{i=1}^{2}N_{(p_{i},q)}\bigcup_{i=1}^{2}N_{(p_{i},2q)}\bigcup N_{(p_{1}p_{2},u)}\notag\\
		 &\quad \bigcup N_{(p_{1}p_{2},2)}\bigcup N_{(p_{1}p_{2},q)}\bigcup N_{(p_{1}p_{2},2q)}\bigcup N_{(u,2)}\bigcup N_{(u,q)}\bigcup N_{(u,2q)}
  \end{align}
For $\bar{x},\bar{y}\in \mathbb{Z}_{p_{1}^{n_1}p_{2}^{n_2}}\oplus\mathbb{Z}_{2^{m_1}q^{m_{2}}}$, consider  following cases:\\
{\bf Case I:} If $\displaystyle \bar{x},\bar{y}\in U(\mathbb{Z}_{p_{1}^{n_1}p_{2}^{n_2}}\oplus\mathbb{Z}_{2^{m_1}q^{m_{2}}})$, then $\bar{x}$ and $\bar{y}$ are adjacent to $\bar{0}$, this gives $d(\bar{x},\bar{y})\leq 2$\\
{\bf Case II:} If $\displaystyle \bar{x}\in U(\mathbb{Z}_{p_{1}^{n_1}p_{2}^{n_2}}\oplus\mathbb{Z}_{2^{m_1}q^{m_{2}}})$ and $\displaystyle \bar{y}\in N_{U}(\mathbb{Z}_{p_{1}^{n_1}p_{2}^{n_2}}\oplus\mathbb{Z}_{2^{m_1}q^{m_{2}}})$, then $\bar{x}=(x_{1},x_{2})$ and for $\bar{y}$, we have following possibilities:\\
		(a) If $\bar{y}\in N_{(p_{i},u)}$, then  $\bar{y}=(\alpha_{1}p_{i},u_{1})$. If $p_{j}\mid x_{1}+\alpha_{1}p_{i}$, for $i\neq j$, then take $\bar{z}=(\alpha_{2}p_{j},0)$, where  $\alpha_{2}p_{j}=p_{i}\alpha_{1}+x_{1}$. Clearly, $\bar{z}$ is adjacent to both $\bar{x}$ and $\bar{y}$.\\
		(b) If $\displaystyle \bar{y}\in N_{(p_{i},2)}$, then $\bar{y}=(\alpha_{1}p_{i},2\beta_{1})$. Take $\bar{z}=(\alpha_{1}p_{i}+\alpha_{2}p_{j},2\beta_{1}+\beta_{2}q)\in U(\mathbb{Z}_{p_{1}^{n_1}p_{2}^{n_2}}\oplus\mathbb{Z}_{2^{m_1}q^{m_{2}}})$, where $i\neq j$ and $p_{i}\nmid\alpha_{2},
		\ 2\nmid\beta_{2}$. Clearly, $\bar{y}$ is adjacent to $\bar{z}$ and $\bar{x}, \bar{z}$ are adjacent to $\bar{0}$. Hence, $d(\bar{x},\bar{y})\leq 3$.\\
		(c) If $\displaystyle \bar{y}\in N_{(p_{i},q)}$, then $\bar{y}=(\alpha_{1}p_{i},\beta_{1}q)$. If $p_{j}\mid x_{1}+\alpha_{1}p_{i}$, then $\alpha_{2}p_{j}= x_{1}+\alpha_{1}p_{i}$. Note that $\bar{y}$ is adjacent to both $\bar{z}_{1}=(\alpha_{2}p_{j},2\beta_{2}), \bar{z}_{2}=(\alpha_{2}p_{j},-2\beta_{2})$ and $\bar{x}$ is either adjacent to $\bar{z}_{1}$ or $\bar{z}_{2}$, from this, we get
		$d(x,y)\leq 3$. \\
		(d) If $\displaystyle \bar{y}\in N_{(p_{i},2q)}$, then $\bar{y}=(\alpha_{1}p_{i},2\beta_{1}q)$. If $p_{j}\nmid x_{1}+\alpha_{1}p_{i}$ then $\bar{x}$ and $\bar{y}$ are adjacent. If  $p_{j}\mid x_{1}+\alpha_{1}p_{i}$, then for $\bar{z}_{1}=(x_{1},2\beta_{1}q)\in N_{(u,2q)}$ and $\bar{z}_{2}=(\alpha_{2}p_{j},x_{2})$, where $\alpha_{2}p_{j}=x_{1}+\alpha_{1}p_{i}$, we get, $[\bar{x},\bar{z}_{1}], [\bar{z}_{1},\bar{z}_{2}]$ and $[\bar{z}_{2},\bar{y}]$ are edges in $G(\mathbb{Z}_{p_{1}^{n_1}p_{2}^{n_2}}\oplus\mathbb{Z}_{2^{m_1},q^{m_2}})$.\\
		(e) If $\bar{y}\in N_{(p_{1}p_{2},2)}$, then $\bar{y}=(\alpha_{1}p_{1}p_{2},2\beta_{1})$. Now, for $\bar{z}=(x_{1},2\beta_{1}+\beta_{2}q)$, where $2\nmid\beta_{2}$, then we have, $[\bar{x},\bar{0}], [\bar{0},\bar{z}]$ and $[\bar{z},\bar{y}]$ are edges in $G(\mathbb{Z}_{p_{1}^{n_1}p_{2}^{n_2}}\oplus\mathbb{Z}_{2^{m_1}q_{1}^{m_2}})$.\\
		(f) If $\bar{y}\in N_{(u,q)}$, then $\bar{y}=(u_{1},\beta_{1}q)$. Clearly, $\bar{y}$ is adjacent to any  $\bar{z}=(\alpha_{2}p_{1}p_{2},2\beta_{2})\in N_{(p_{1}p_{2},2)}$. Note that, $\bar{x}$ is adjacent to either $\bar{z}$ or $-\bar{z}$, this implies $d(\bar{x},\bar{y})\leq 3$.\\
		With similar techniques, we can prove that for all other possibilities $d(\bar{x},\bar{y})\leq 3$.\\ 
		{\bf Case III:} Let $\displaystyle \bar{x},\bar{y}\in N_{U}(\mathbb{Z}_{p_{1}^{n_2}p_{2}^{n_2}}\oplus\mathbb{Z}_{2^{m_1}q^{m_2}})$. We prove for some possibilities for $\bar{x}$ and $\bar{y}$, we give proof for some possibilities and all other possibilities follows from the same procedure:\\
		(a) If $\bar{x},\bar{y}\in N_{(p_i,u)}$, then for any $\bar{z}=(\alpha_{3}p_{j},0)\in N_{(p_{j},2q)}$, for $i\neq j$, we have $\bar{x}$ and $\bar{y}$ are adjacent to $\bar{z}$. Similarly, if both $\bar{x}$ and $\bar{y}$ in some set in (\ref{eq2}), then there exist $\bar{z}$ in any other set in (\ref{eq2}) such that $[\bar{x},\bar{z}]$ and $[\bar{z},\bar{y}]$ are edges in $\displaystyle G(\mathbb{Z}_{p_{1}^{n_1}p_{2}^{n_2}}\oplus\mathbb{Z}_{2^{m_1}q^{m_{2}}})$. \\
		(b) If $\bar{x}\in N_{(p_{i},u)}$ and $\bar{y}\in N_{(p_{j},2)}$, then $\bar{x}=(\alpha_{1}p_{i},u_{1})$ and $\bar{y}=(\alpha_{2}p_{j},2\beta_{1})$. If $i\neq j$, then take an element $\bar{z}_{1}=(\alpha_{1}p_{i}+\alpha_{2}p_{j},2\beta_{1}+\beta_{2}q)$, where $2\nmid \beta_{2}$. From this, we have, $\bar{y}$ is  adjacent to $\bar{z}_{1}$ and $\bar{z}_{1},\bar{x}$ are adjacent to $\bar{z}_{2}=(\alpha_{1}p_{i}+\alpha_{2}p_{j},0)$. If $i=j$ then take an element $\bar{z}_{1}=(\alpha_{3}p_{k},2\beta_{1}+\beta_{2}q)\in N_{(p_{k},u)}$, where $2\nmid \beta_{2}$ and $k\neq i,j$. From this, we have, $\bar{y}$ is  adjacent to $\bar{z}_{1}$ and $\bar{z}_{1},\bar{x}$ are adjacent to $\bar{z}_{2}=(\alpha_{1}p_{i}+\alpha_{3}p_{k},0)$. \\
		(c) If $\bar{x}\in N_{(p_{i},u)}$ and $\bar{y}\in N_{(p_{j},q)}$, then $\bar{x}=(\alpha_{1}p_{i},u_{1})$ and $\bar{y}=(\alpha_{2}p_{j},\beta_{1} q)$. If $i=j$, then for any element $\bar{z}=(\alpha_{3}p_{k},2\beta_{2})\in N_{(p_{k},2)}$, where $k\neq i,j$. Clearly, $\bar{y}$ adjacent to $\bar{z}$. Note that $\bar{x}$ is either adjacent to $\bar{z}$ or $-\bar{z}$. If $i\neq j$, then for elements $\bar{z}_{1}=(\alpha_{1}p_{i},2\beta_{2}), \bar{z}_{2}=(\alpha_{1}p_{i},-2\beta_{2})\in N_{(p_{i},2)}$, we get, $\bar{y}$ is adjacent to both $\bar{z}_{1} \ \& \bar{z}_{2}$ and $\bar{x}$ is either adjacent to $\bar{z}_{1}$ or $\bar{z}_{2}$ \\
	(d) If $\bar{x}\in N_{(p_{i},u)}$ and $\bar{y}\in N_{(p_{1}p_{2},2q)}$ then $\bar{x}=(\alpha_{1}p_{i},u_{1})$ and $\bar{y}=(\alpha_{2}p_{1}p_{2},2\beta_{1}q)$. Take $\bar{z}_{1}=(\alpha_{1}p_{i}+\alpha_{3}p_{j},2\beta_{1}q)$ is adjacent to $\bar{x}$ and $\bar{z}_{2}=(\alpha_{1}p_{i}+\alpha_{3}p_{j},u_{1})$ is adjacent to both $\bar{z}_{1}$ and $\bar{y}$.\\
	(e) If $\bar{x}\in N_{(p_{i},u)}$ and $\bar{y}\in N_{(u,2)}$, then $\bar{x}=(\alpha_{1}p_{i},u_{1})$ and $\bar{y}=(u_{2},2\beta_{1})$. Consider, the case $p_{j}\mid u_{2}+\alpha_{1}p_{i}$, which implies $\alpha_{2}p_{j}= u_{2}+\alpha_{1}p_{i}$. For this, take an element $\bar{z}_{1}=(\alpha_{2}p_{j},2\beta_{1}+\beta_{2}q)$, where $2\nmid \beta_{2}$. Clearly $\bar{y}$ is adjacent to $\bar{z}_{1}$ and $\bar{z}_{1},\bar{x}$ are adjacent to $\bar{z}_{2}=(\alpha_{2}p_{j},0)$.\\
	(f) If $\bar{x}\in N_{(p_{i},2)}$ and $\bar{y}\in N_{(p_{j},2q)}$, then $\bar{x}=(\alpha_{1}p_{i},2\beta_{1})$ and $\bar{y}=(\alpha_{2}p_{j},2\beta_{2}q)$. If $i=j$, then take, $\bar{z}=(\alpha_{3}p_{k},2\beta_{1}+\beta_{3}q)$, where $k\neq i,j$ and $2\nmid \beta_{3}$. Note that $\bar{z}$ is adjacent to both $\bar{x}$ and $\bar{y}$. For $i\neq j$, take $\bar{z}=(\alpha_{1}p_{i}+\alpha_{2}p_{j},2\beta_{1}+\beta_{3}q)$.\\
	(g) If $\bar{x}\in N_{(p_{i},q)}$ and $\bar{y}\in N_{(p_{j},2q)}$, then $\bar{x}=(\alpha_{1}p_{i},\beta_{1}q)$ and $\bar{y}=(\alpha_{2}p_{j},2\beta_{2}q)$. If $i\neq j$, then $\bar{z}_{1}=(\alpha_{1}p_{i}+\alpha_{2}p_{j},2\beta_{3})$, where $q\nmid \beta_{3}$ is adjacent to $\bar{x}$ and $\bar{z}_{2}=(\alpha_{1}p_{i}+\alpha_{2}p_{j},\beta_{1}q+2\beta_{3})$ is adjacent to both $\bar{z}_{1}$ and $\bar{y}$. Other case also follows the same procedure as above.\\
	(h) If $\bar{x}\in N_{(p_{i},q)}$ and $\bar{y}\in N_{(p_{1}p_{2},2q)}$, then $\bar{x}=(\alpha_{1}p_{i},\beta_{1}q)$ and $\bar{y}=(\alpha_{2}p_{1}p_{2},2\beta_{2}q)$. Take $\bar{z}_{1}=(\alpha_{1}p_{i}+\alpha_{3}p_{j},2\beta_{3})$ is adjacent to $\bar{x}$ and $\bar{z}_{2}=(\alpha_{1}p_{i}+\alpha_{3}p_{j},\beta_{1}q+2\beta_{3})$ is adjacent to both $\bar{z}_{1}$ and $\bar{y}$, where $p_{j}\nmid \alpha_{3}, \ q\nmid \beta_{3}$.\\
	 (i) If $\bar{x}\in N_{(p_{i},2q)}$ and $\bar{y}\in N_{(p_{1}p_{2},u)}$, then $\bar{x}=(\alpha_{1}p_{i},2\beta_{1}q)$ and $\bar{y}=(\alpha_{2}p_{1}p_{2},u_{1})$. Consider an elements, $\bar{z}_{1}=(\alpha_{1}p_{i}+\alpha_{3}p_{j},u_{1})$ and $\bar{z}_{2}=(\alpha_{1}p_{i}+\alpha_{3}p_{j},2\beta_{1}q)$, where $p_{i}\nmid \alpha_{3}$. We get $\bar{z}_{1}$ is adjacent to $\bar{x}$ and $\bar{z}_{2}$ is adjacent to both $\bar{z}_{1}$ and $\bar{y}$.\\
	  (j) If $\bar{x}\in N_{(p_{1}p_{2},2)}$ and $\bar{y}\in N_{(p_{1}p_{2},q)}$, then $\bar{x}=(\alpha_{1}p_{1}p_{2},2\beta_{1})$ and $\bar{y}=(\alpha_{2}p_{1}p_{2},\beta_{2}q)$. Consider an element $\bar{z}_{1}=(u,\beta_{2}q)\in N_{(u,q)}$ adjacent to $\bar{x}$ and $\bar{z}_{2}=(u,2\beta_{1})$ is adjacent to both $\bar{z}_{1}$ and $\bar{y}$.\\
	  (k) If $\bar{x}\in N_{(u,2)}$ and $\bar{y}\in N_{(u,q)}$,  then $\bar{x}=(u_{1},2\beta_{1})$ and $\bar{y}=(u_{2},\beta_{2}q)$. Consider an element $\bar{z}_{1}=(\alpha_{1}p_{1}p_{2},\beta_{2}q)\in N_{(p_{1}p_{2},q)}$ adjacent to $\bar{x}$. Also, consider $\bar{z}_{2}=(u_{2},2\beta_{1})\in N_{(u,2)}$, which is adjacent to both $\bar{z}_{1}$ and $\bar{y}$.\\
	  Hence, $d(\bar{x}, \bar{y})\leq 3$, for all $\bar{x},\bar{y}\in \mathbb{Z}_{p_{1}^{n_1}p_{2}^{n_2}} \oplus \mathbb{Z}_{2^{m_1}q^{m_2}}$ and $G(\mathbb{Z}_{p_{1}^{n_1}p_{2}^{n_2}} \oplus \mathbb{Z}_{2^{m_1}q^{m_2}})$ is a connected graph, which
	  gives diam$(G(\mathbb{Z}_{p_{1}^{n_1}p_{2}^{n_2}} \oplus \mathbb{Z}_{2^{m_1}q^{m_2}}))\leq 3$.
	\end{proof}
\end{theorem}
\begin{corollary}\label{thm38}
	Let $\displaystyle G(\mathbb{Z}_{p_{1}^{n_1}p_{2}^{n_2}} \oplus \mathbb{Z}_{2^{m_1}q^{m_2}})$ be a unit graph, where $\displaystyle p_{1},p_{2}$ and $q$ are odd primes. Then $\displaystyle \lambda(G(\mathbb{Z}_{p_{1}^{n_1}p_{2}^{n_2}} \oplus \mathbb{Z}_{2^{m_1}q^{m_2}}))=2^{m_{1}-1}\phi(p_{1}^{n_1}p_{2}^{n_2})\phi(q^{m_2})$.
\begin{proof}
It follows from, Theorem \ref{thm4} and \ref{thm34}.
\end{proof}
\end{corollary}

\begin{theorem} Let $\displaystyle G(\mathbb{Z}_{n}\oplus\mathbb{Z}_{m})$ be a unit graph, where $n=p^{n_{1}}_{1}p^{n_{2}}_{2}, m=q_{1}^{m_1}q_{2}^{m_2}$ and $p_{1},\ p_{2} \ q_{1} \ \& \ q_{2}$ are  primes. Let $H$ be a $|V|\times |E|$ incidence matrix of $\displaystyle G(\mathbb{Z}_{n}\oplus\mathbb{Z}_{m})$.
    \begin{enumerate}
      \item If both $\displaystyle m$ and $n$ are odd, then  $\displaystyle C_{2}(H)=\left[\frac{(mn-1)\phi(m)\phi(n)}{2},mn-1,\phi(m)\phi(n)-1\right]_{2}$ is the  binary code generated by $H$ over finite field $\mathbb{F}_{2}$.
      \item If $m$ is even and $n$ is odd, then for any odd prime $r$, \\ $\displaystyle C_{r}(H)=\left[\frac{mn\phi(m)\phi(n)}{2},mn-1,\phi(m)\phi(n)\right]_{r}$ is  the $r$-ary code generated by $H$ over finite field $\mathbb{F}_{r}$.
    \end{enumerate}
\begin{proof}
\begin{enumerate}
  \item If $m$ and $n$ both are odd, then $p_{1},p_{2},q_{1}$ and $q_{2}$ are odd primes. By Theorem (\ref{thm31}), $\displaystyle G(\mathbb{Z}_{n}\oplus\mathbb{Z}_{m})$ is a connected graph and hence by Theorem (\ref{thm2}), binary code generated by $H$ is
      $\displaystyle C_{2}(H)=[|E|,|V|-1,\lambda(G(\mathbb{Z}_{n}\oplus\mathbb{Z}_{m}))]_{2}$. Now from Theorem (\ref{thm41}) and Corollaey (\ref{thm33}), we get $\displaystyle |E|= \frac{(mn-1)\phi(m)\phi(n)}{2},\\  |V|-1=mn-1$ and $\displaystyle \lambda(G(\mathbb{Z}_{n}\oplus \mathbb{Z}_{m}))=\phi(m)\phi(n)-1$.
  \item If $m$ is even and $n$ is odd, then either $q_{1}=2$ or $q_{2}=2$. By Theorem (\ref{thm34}) and Lemma (\ref{thm24}), $\displaystyle G(\mathbb{Z}_{n}\oplus \mathbb{Z}_{n})$ is a connected bipartite graph and hence by Theorem (\ref{thm5}), $r$-ary code generated by $H$ is $\displaystyle C_{r}(H)=[|E|,|V|-1,\lambda(G(\mathbb{Z}_{n}\oplus \mathbb{Z}_{m}))]_{r}$. Now using Theorem (\ref{thm41}) and Corollary (\ref{thm38}), we conclude the result.
\end{enumerate}
\end{proof}
\end{theorem}

\begin{corollary}
Let $C_{r}(H)$ and $C_{2}(H)$ denote the linear codes generated from incidence matrices of $\displaystyle G(\mathbb{Z}_{p_{1}^{n_1}p_{2}^{n_2}}\oplus\mathbb{Z}_{2^{m_1}q^{m_{1}}})$ and $\displaystyle G(\mathbb{Z}_{p_{1}^{n_1}p_{2}^{n_2}}\oplus\mathbb{Z}_{q_{1}^{m_1}q^{m_{1}}})$. Then
\begin{enumerate}
  \item Dual of code $C_{2}$ is $\displaystyle C^{\perp}_{2}=\left[\frac{(mn-1)\phi(m)\phi(n)}{2},\frac{(mn-1)[\phi(m)\phi(n)-2]}{2},3\right]_{2}$, where $\displaystyle n=p^{n_{1}}_{1}p^{n_{2}}_{2}$ and $\displaystyle m=q^{m_{1}}_{1}q^{m_{2}}_{2}$.
  \item Dual of code $C_{r}$ is $\displaystyle C^{\perp}_{r}=\left[\frac{mn\phi(m)\phi(n)}{2},\frac{mn(\phi(m)\phi(n)-2)+2}{2},4\right]_{r}$, where $\displaystyle m=2^{m_1}q^{m_{2}}_{1}$ and $n=p_{1}^{n_1}p_{2}^{n_2}$.
\end{enumerate}
\begin{proof}
\begin{enumerate}
\item From Theorem (\ref{thm8}), $\text{dim}(\displaystyle C^{\perp}_{2})=\frac{(mn-1)[\phi(m)\phi(n)-2]}{2}$ and from Theorem \ref{thm6}, we have $\displaystyle d(C^{\perp}_{2})=g_{r}(G(\mathbb{Z}_{p_{1}^{n_1}p_{2}^{n_2}}\oplus \mathbb{Z}_{q_{1}^{m_1}q_{2}^{m_2}}))$. Since, $(p_{1}+p_{2},q_{1}+q_{2})$ is adjacent to both $  (p_{1},q_{1}), \ (p_{2},q_{2})$ and $(p_{1},q_{1})$ is adjacent to  $(p_{2},q_{2})$,  we conclude the result.
  \item From Theorem (\ref{thm8}), $\text{dim}(\displaystyle C^{\perp}_{2})=\frac{mn(\phi(m)\phi(n)-2)+2}{2}$ and from Theorem \ref{thm6}, we have $\displaystyle d(C^{\perp}_{2})=g_{r}(G(\mathbb{Z}_{p_{1}^{n_1}p_{2}^{n_2}}\oplus \mathbb{Z}_{2^{m_1}q_{2}^{m_2}}))$. From Theorem (\ref{thm24}) $\displaystyle G(\mathbb{Z}_{p_{1}^{n_1}p_{2}^{n_2}}\oplus\mathbb{Z}_{2^{m_1}q^{m_{2}}})$ is a bipartite graph and hence it has grith as even number. Consider an elements $\bar{z}_{1}=(\alpha_{1}p_{1},2\beta_{1})\in N_{(p_{1},2)}, \bar{z}_{2}=(\alpha_{2}p_{2},\beta_{1}q)\in N_{(p_{2},q)} \ \bar{z}_{3}=(\alpha_{1}p_{1}+\alpha_{2}p_{2},2\beta_{2})\in N_{(u,2)}$ and $\bar{z}_{4}=(\alpha_{1}p_{1}+\alpha_{2}p_{2},\beta_{1}q)\in N_{(u,q)}$. From this we have, $[\bar{z}_{1},\bar{z}_{2}], \ [\bar{z}_{2},\bar{z}_{3}][\bar{z}_{3},\bar{z}_{4}] $ and $[\bar{z}_{4},\bar{z}_{1}]$ are edges in $G(\mathbb{Z}_{p_{1}^{n_1}p_{2}^{n_2}}\oplus\mathbb{Z}_{2^{m_1}q^{m_{2}}})$.
\end{enumerate}
\end{proof}
\end{corollary}
Based on the procedure to obtain the results  in Section 4 \& 5, we state following conjectures\\
{\bf Conjecture I:}
Let $G(\mathbb{Z}_{n}\oplus\mathbb{Z}_{m})$ be a unit graph.
\begin{enumerate}
	\item If both $m$ and $n$ are odd, then  $G(\mathbb{Z}_{n}\oplus\mathbb{Z}_{m})$ is a connected graph and $\text{diam}(G(\mathbb{Z}_{n}\oplus\mathbb{Z}_{m}))\leq 2$.
	\item If exactly one of $m$ and $n$ is even, then  $G(\mathbb{Z}_{n}\oplus\mathbb{Z}_{m})$ is a connected graph and $\text{diam}(G(\mathbb{Z}_{n}\oplus\mathbb{Z}_{m}))\leq 3$.
\end{enumerate}
{\bf Conjecture II:}
Let $\displaystyle G(\mathbb{Z}_{n}\oplus\mathbb{Z}_{m})$ be a unit graph and $H$ be a $|V|\times |E|$ incidence matrix of $G(\mathbb{Z}_{n}\oplus\mathbb{Z}_{m})$.
\begin{enumerate}
	\item If both $m$ and $n$ are odd, then  $\displaystyle C_{2}(H)=\left[\frac{(mn-1)\phi(m)\phi(n)}{2},mn-1,\phi(m)\phi(n)-1\right]_{2}$ is the binary code generated by $H$ over the finite field $\mathbb{F}_{2}$.
	\item If exactly one from $m$ and $n$ is even, then for any odd prime $r$, $\displaystyle C_{r}(H)=\left[\frac{mn\phi(m)\phi(n)}{2},mn-1,\phi(m)\phi(n)\right]_{r}$ is the $r$-ary code generated by $H$ over the finite field $\mathbb{F}_{r}$.
\end{enumerate}
\section{Conclusion}
In this paper, we constructed $r$-ary linear codes from the  incidence matrices of unit graphs $\displaystyle G(\mathbb{Z}_{n}\oplus \mathbb{Z}_{m})$, where $n$ and $m$ being the power of primes or product of powers of two primes. Furthermore, we found minimum distance of corresponding dual codes over finite field $\displaystyle \mathbb{F}_{r}$. We state two conjectures on construction of linear codes from unit graphs $\displaystyle G(\mathbb{Z}_{n}\oplus \mathbb{Z}_{m})$ for any $m$ and $n$. Examine the permutation decoding techniques, covering radius of constructed codes and one can construct linear codes from unit graph over different commutative rings is the further scope to work.

\end{document}